\documentclass[final,leqno]{siamltex}

\usepackage{epsfig,amsmath,amssymb,latexsym,color}
\usepackage{subfigure,colordvi}
\usepackage[curve]{xypic}

\def\RR{\mathbb{R}}

\def\bfe{{\bf e}}
\def\bff{{\bf f}}

\def\bfn{{\bf n}}
\def\bfp{{\bf p}}
\def\bfR{{\bf R}}

\def\bfu{{\bf u}}
\def\bfv{{\bf v}}
\def\bfw{{\bf w}}

\def\bfx{{\bf x}}
\def\bfy{{\bf y}}

\def\det{{\textsf{det}}}
\def\eps{{\epsilon}}

\newcommand{\eref}[1]{(\ref{#1})}

\newcommand{\tref}[1]{Table~\ref{#1}}
\newcommand{\vm}[1]{\mathbf{#1}}
\newcommand{\vx}{\vm{x}}

\newcommand{\cancel}[1]{ }

\newcommand{\raw}{\rightarrow}
\newcommand{\diam}{\textnormal{diam}}
\newcommand{\vn}[1]{\left\|#1\right\|} 
\newcommand{\vsn}[1]{\left|#1\right|} 
\newcommand{\ds}{\displaystyle}

\newcommand{\overbar}[1]{\mkern 1.5mu\overline{\mkern-1.5mu#1\mkern-1.5mu}\mkern 1.5mu}

\begin{document}

\title{Gradient bounds for Wachspress coordinates on polytopes}
\author{Michael S. Floater\thanks{
   Department of Mathematics,
    University of Oslo,
    P.O.\ Box 1053 Blindern, 0316 Oslo, Norway,
   {\tt michaelf@ifi.uio.no}} \and
  Andrew Gillette\thanks{Department of Mathematics, University of California, San Diego, 9500 Gilman Drive MC 0112, La Jolla, CA 92093 {\tt akgillette@mail.ucsd.edu}}
  \and
   N. Sukumar\thanks{Department of Civil \& Environmental Engineering, University of California, Davis, One Shields Avenue, Davis, CA 95616 {\tt nsukumar@ucdavis.edu}}
}


\maketitle

\begin{abstract}
We derive upper and lower bounds on the gradients of Wachspress coordinates defined over any simple convex $d$-dimensional polytope $P$.
The bounds are in terms of a single geometric quantity $h_\ast$, which denotes the minimum distance between a vertex of $P$ and any hyperplane containing a non-incident face.
We prove that the upper bound is sharp for $d=2$ and analyze the bounds in the special cases of hypercubes and simplices.
Additionally, we provide an implementation of the Wachspress coordinates on convex polyhedra using Matlab and employ them in a 3D finite element solution of the 
Poisson equation on a non-trivial polyhedral mesh.
As expected from the upper bound derivation, the $H^1$-norm of the error in the method converges at a linear rate with respect to the size of the mesh elements.
\end{abstract}



\pagestyle{myheadings}
\thispagestyle{plain}

\section{Introduction}\label{sec:intro}

Given a set of generalized barycentric coordinates $\{\phi_\bfv\}$ on a polytope $P$, viewed as an open set in $\RR^d$ with vertex set $V$, the standard vertex-based interpolation $I$ of a function $u:\overbar P\raw\RR$ is given by
\begin{equation}
\label{eq:I-def}
I(u) := \sum_{\bfv\in V} u(\bfv)\phi_\bfv.
\end{equation}
Let $|\cdot|$ denote the standard Euclidean norm.
Observe that if $\phi_\bfv \in C^1(P)$ for all $\bfv\in V$, then
\[ \sup_{\bfx\in P}|\nabla I(u)(\bfx)|  
  \leq \sup_{\bfx\in P}\sum_{\bfv\in V} \left|u(\bfv) \nabla\phi_\bfv(\bfx)\right| 
  \leq \Lambda\;\max_{\bfv\in V}|u(\bfv )|, \]
where
\begin{equation}
\label{eq:Lambda-def}
\Lambda := \sup_{\bfx \in P} \lambda(\bfx)\quad\text{with}\quad\lambda(\bfx) := \sum_{\bfv\in V} |\nabla \phi_\bfv(\bfx)|.
\end{equation}
The main result of this paper is an upper bound on $\Lambda$ over the class of simple convex polytopes when the $\phi_\bfv$ coordinates in $I$ are generalized Wachspress coordinates.
We also derive lower bounds on $\Lambda$ to illustrate the sharpness of the upper bound and provide code written in Matlab for numerical experimentation with the Wachspress coordinates on polygons and polyhedra.

Our motivation for this analysis stems from the growing interest in using  generalized barycentric coordinates for finite element methods on polygonal and polyhedral meshes~\cite{dBCMMR13,GRB1,RGB1,ST2004,WBG}.
In such methods, the interpolant $I$ is suitable if it admits an \textit{a priori} error estimate of the form
\begin{equation}
\label{eq:opt-conv-est}
\vn{u- I(u)}_{H^1(P)} \leq C~\diam(P)\vsn{u}_{H^2(P)}\qquad\forall u\in H^2(P),
\end{equation}
where $H^k(P)$ denotes the degree $k$ Sobolev space over $P$.
Following~\cite[Section 4]{GRB1} and classical finite element sources~\cite{BS08,Strang1973}, we can take the constant $C$ to be 
\begin{equation}
\label{eq:C-def}
C := (1+C_S(1+\Lambda))\sqrt{1+C_{BH}^2},
\end{equation}
with $C_S$ the Sobolev embedding constant satisfying $\vn{u}_{C^0(\overline{P})}\leq C_S\vn{u}_{H^2(P)}$ independent of $u\in H^2(P)$ and $C_{BH}$ the Bramble-Hilbert constant for linear approximation on a class of polytopes of diameter 1.
Therefore, to prove (\ref{eq:opt-conv-est}) in this context, it suffices to provide an upper bound on $\Lambda$ holding over the class of polytopes to be used as domain mesh elements and the set of coordinate functions $\phi_\bfv$ to be used for interpolation via $I$.
Here we consider the class of polytopes of dimension $d$ that are \textit{simple}, meaning the number of faces incident to each vertex is exactly $d$.

We summarize our results in Table~\ref{tab:bounds}.
Our bounds are in terms of a single geometric quantity $h_\ast$, which denotes the minimum distance between a vertex of $P$ and any hyperplane containing a non-incident face.
The bounds are inversely proportional to $h_\ast$, representing the fact that geometries with small $h_\ast$ values can result in large values of constant C in the \textit{a priori} error estimate (\ref{eq:opt-conv-est}).
The usefulness of $h_\ast$ as a measure of geometric quality is discussed further in Section~\ref{sec:conc}.

\begin{table}[ht]
\centering
\sbox{\strutbox}{\rule{0pt}{0pt}}           
\begin{tabular}[.75\textwidth]{@{\extracolsep{\fill}} rcccccc}
simple convex polytope in $\RR^d$ && $\ds \frac 1{h_\ast}$ & $\leq$ & $\Lambda$ & $\leq$ & $\ds \frac{2d}{h_\ast}$ \\
\\[3mm]
\hline \\[3mm]
$d$-simplex in $\RR^d$ && $\ds \frac 1{h_\ast}$ & $\leq$ & $\Lambda$ & $\leq$ & $\ds \frac{d+1}{h_\ast}$ \\
\\[3mm]
\hline \\[3mm]
hyper-rectangle in $\RR^d$ && $\ds \frac 1{h_\ast}$ & $\leq$ & $\Lambda$ & $\leq$ & $\ds\frac{d+\sqrt d}{h_\ast}$ \\
\\[3mm]
\hline \\[3mm]
regular $n$-gon in $\RR^2$ && $\ds \frac {2(1+\cos(\pi/n))}{h_\ast}$ & $\leq$ & $\Lambda$ & $\leq$ & $\ds \frac{4}{h_\ast}$\\
\\[3mm]
\hline \\[3mm]
\end{tabular}
\caption{\noindent Our bounds on $\Lambda$ for various polytope types are summarized above.  In the case of simplices and hyper-rectangles, $\Lambda$ obtains the upper bound when the shapes are regular.  The lower bound for regular $n$-gons approaches the upper bound as $n\raw\infty$, meaning $2d/h_\ast$ is a sharp upper bound in the case $d=2$.}
\label{tab:bounds}
\end{table}

The lower bounds on $\Lambda$ hold for any generalized barycentric coordinates $\{\phi_\bfv\}$ that are $C^1$ at the vertices of $P$, while the upper bounds hold in the specific case of generalized Wachspress coordinates.
The namesake of the Wachspress coordinates is the author of the book where they first appeared in the context of polygonal finite elements~\cite{Wach}; see also \cite{Wach2}.
Warren first generalized this definition to polytopes in~\cite{Warr} and later to convex sets with coauthors in~\cite{WSHD}.
Our notation in this work follows the definition of the coordinates
in Warren et al.~\cite{WSHD} and Ju et al.~\cite{JSWD}.

The outline of the paper is as follows.
In Section~\ref{sec:upper-bds}, we define generalized barycentric coordinates, generalized Wachspress coordinates, and $h_\ast$ precisely before deriving the upper bounds given in Table~\ref{tab:bounds}.
In Section~\ref{sec:lower-bds}, we derive the lower bounds given in Table~\ref{tab:bounds}.
In Section~\ref{sec:results}, we present an implementation of the Wachspress coordinates on convex polyhedra and employ them in a 3D finite element solution of the Poisson equation on a non-trivial polyhedral mesh.
We present our conclusions and discuss future directions in Section~\ref{sec:conc}.

\section{Upper bounds}\label{sec:upper-bds}

\begin{figure}[h]
\centering
\includegraphics[width=.5\textwidth]{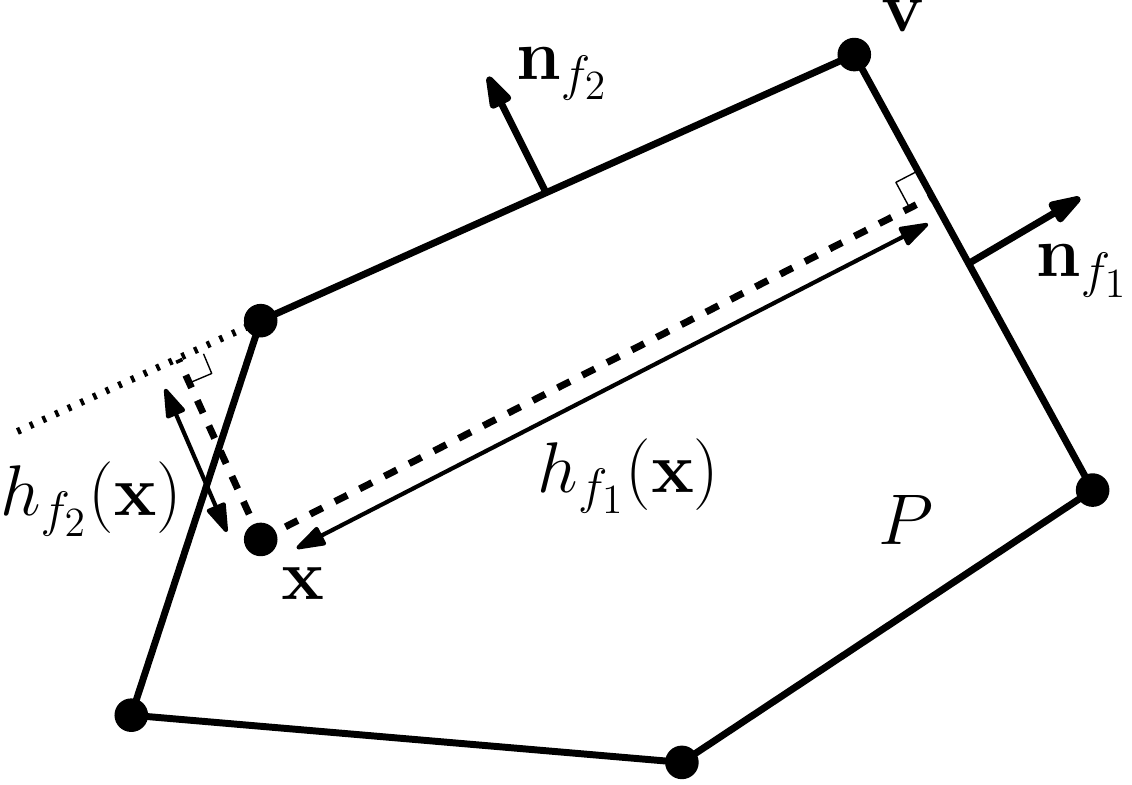}
\noindent \caption{The value of the generalized Wachspress basis function $\phi_\bfv$ at a point $\bfx$ inside a convex polytope $P$ is defined in terms of the location of $\bfv$, the normals $\bfn_{f_i}$ to the faces $f_i$ incident to $\bfv$, and the distances $h_i(\bfx)$ to the planes containing the $f_i$.
}
\label{fg:notation-2d}
\end{figure}

We start by fixing notation for describing polytope geometry.
Let $P \subset \RR^d$ be a convex $d$-dimensional polytope, viewed as an open set, with $V$ and $F$ the sets of vertices and $(d-1)$-dimensional faces of $\partial P$, respectively. 
Let $V_f \subset V$ denote the set of vertices of face $f \in F$, and $F_\bfv \subset F$ denote the set of faces incident to the vertex $\bfv \in V$.
Assume that $P$ is \textit{simple} meaning $|F_\bfv|=d$ for all $\bfv\in V$.

For any $\bfx$ in $P$, let $h_f(\bfx)$ denote the perpendicular distance from $\bfx$ to the $(d-1)$-hyperplane in $\RR^d$ containing the face $f\in F$.
Letting $\bfn_f$ denote the outward unit normal to $f$ and $\bfv$ any vertex in $f$, we can compute $h_f(\bfx)$ via
\begin{equation}
\label{eq:def-hf}
h_f(\bfx) = (\bfv - \bfx) \cdot \bfn_f.
\end{equation}
We will also make use of scaled normal vectors defined by
\begin{equation}
\label{eq:def-pf}
\bfp_f(\bfx) := \bfn_f / h_f(\bfx).
\end{equation}
The \textit{generalized Wachspress coordinates} for a simple polytope $P$ as above are the functions $\phi_\bfv : P \to \RR$, for $\bfv \in V$, given by the formula,
\begin{equation}
\label{eq:lambda-def}
\phi_\bfv(\bfx) := \frac{w_\bfv(\bfx)}{ W(\bfx)},
\end{equation}
where
\begin{equation}
\label{eq:wachsimple}
W(\bfx) := \sum_{\bfu \in V} w_\bfu(\bfx),
\quad\text{and}\quad
 w_\bfv(\bfx) := \det(\bfp_{f_1}(\bfx), \cdots, \bfp_{f_d}(\bfx)),
\end{equation}
where $f_1,\ldots,f_d$ are the $d$ faces adjacent to $\bfv$  listed in an counter-clockwise ordering around $\bfv$ as seen from outside $P$, and $\det$ denotes the regular vector determinant in $\RR^d$.
The notation is summarized for the case $d=2$ in Figure~\ref{fg:notation-2d}.


We have $h_\bfv(\bfx)>0$ on $P$ and $\det(\bfn_{f_1},\cdots,\bfn_{f_d})>0$ by the strict convexity assumption and incident face ordering convention.
Thus, $w_{\bfv}(\bfx)>0$, $W(\bfx)>0$, and hence $\phi_{\bfv}(\bfx)>0$ on $P$.
The partition of unity property $\sum_{\bfv\in V} \phi_{\bfv}(\bfx) = 1$ is immediate from (\ref{eq:lambda-def}) and a proof of the linear precision property $\sum_{\bfv\in V} \phi_\bfv(\bfx) \bfv = \bfx$ can be found in \cite{WSHD} and \cite{JSWD}.
The \textit{linear completeness property} then follows immediately: for any linear function $L:P\rightarrow\RR$,
\begin{equation}
\label{eq:linprec3d}
\sum_{\bfv \in V} \phi_\bfv(\bfx) L(\bfv) = L(\bfx).
\end{equation}

Before bounding $\Lambda$, we first derive a convenient expression for 
$\nabla \phi_\bfv$ in terms of the coordinates $\phi_\bfu$
and the (vector-valued) ratios
\begin{equation}
\label{eq:def-R}
\bfR_\bfv(\bfx) := \frac{\nabla w_\bfv(\bfx)}{w_\bfv(\bfx)}.
\end{equation}

\begin{lemma}\label{lem:R}
For $\bfv \in V$,
\begin{equation}\label{eq:R3d}
 \nabla \phi_\bfv = \phi_\bfv
                      (\bfR_\bfv - \sum_{\bfu\in V} \phi_\bfu \bfR_\bfu).
\end{equation}
\end{lemma}\\

\emph{Proof.}
Taking the gradient of $\phi_\bfv$ yields
\[ \nabla \phi_\bfv 
  = \frac{\nabla w_\bfv}{W} - \frac{w_\bfv \nabla W}{W^2}
  = \phi_\bfv \bfR_\bfv  - \phi_\bfv \frac{\nabla W}{W}. \]
The result follows from the observation that
\[ \frac{\nabla W}{W}
         = \sum_{\bfu\in V} \frac{\nabla w_\bfu}{W}
         = \sum_{\bfu\in V} \phi_\bfu \bfR_\bfu. \qquad \endproof \] 

The bound on $\Lambda$ will be in terms of the minimum distance between a vertex of $P$ and any hyperplane containing a non-incident face.
We denote this geometric quantity by 
\begin{equation}
\label{eq:hstar-def}
 h_\ast := \min_{f \in F} \min_{\bfu \in V \setminus V_f} h_f(\bfu).
\end{equation}
We also introduce the notation $\mu_f$ to denote the sum of the coordinates associated with the face $f$, i.e.\
\begin{equation}
\label{eq:muf-def}
\mu_f(\bfx) := \sum_{\bfv \in V_f} \phi_\bfv(\bfx).
\end{equation}

\begin{lemma}
\label{lem:hstarineq}
For $f \in F$,
\begin{equation}\label{eq:hstarineq3d}
 \frac{1-\mu_f(\bfx)}{h_f(\bfx)} \le \frac{1}{h_\ast}.
\end{equation}
\end{lemma}\\

\emph{Proof.}
Since $h_f: P\raw\RR$ is linear, the linear completeness property (\ref{eq:linprec3d}) implies that
\[ h_f(\bfx) = \sum_{\bfu \in V} \phi_\bfu(\bfx) h_f(\bfu)
  = \sum_{\bfu \in V \setminus V_f} \phi_\bfu(\bfx) h_f(\bfu)
   \ge h_\ast \sum_{\bfu \in V \setminus V_f} \phi_\bfu(\bfx)
\qquad \endproof \]

\begin{theorem}
\label{thm:bound-gen-d}
Let $P$ be a simple convex polytope in $\RR^d$ and let $\Lambda$ from (\ref{eq:Lambda-def}) be defined using generalized Wachspress coordinates from (\ref{eq:lambda-def}).
Then
\begin{equation}\label{eq:bound-gen-d}
 \Lambda \leq \frac{2d}{h_\ast}.
\end{equation}
\end{theorem}\\

\begin{proof}
We first compute the gradient of $w_\bfv$ as defined in (\ref{eq:wachsimple}).
Observing that $\nabla h_f(\bfx) = -\bfn_f$ for $\bfx \in P$,
\begin{equation}
\label{eq:pf-obs1}
 \nabla \left(\frac{1}{\prod_{\ell=1}^d h_{f_\ell}}\right) =
\frac{1}{\prod_{\ell=1}^d h_{f_\ell}}
   \left(\sum_{\ell=1}^d\frac{\bfn_{f_\ell}}{h_{f_\ell}}\right).
\end{equation}
Thus, recalling the definition of $\bfp_{f_i}$ from (\ref{eq:def-pf}), we have 
\begin{equation}
\label{eq:pf-obs2}
 \nabla w_\bfv
  = \nabla\left(\frac{\det(\bfn_{f_1},\cdots,\bfn_{f_d})}{\prod_{\ell=1}^d h_{f_\ell}}\right)
  = w_\bfv \sum_{\ell=1}^d \bfp_{f_\ell}
  = w_\bfv 
       \sum_{f \in F_\bfv} \bfp_f.
\end{equation}
Recalling the definition of $\bfR_\bfv$ from (\ref{eq:def-R}), we have just shown
\begin{equation}
\label{eq:Rww3d}
 \bfR_\bfv = \sum_{f \in F_\bfv} \bfp_f.
\end{equation}
Hence, by Lemma~\ref{lem:R},
\begin{align*}
 \frac{\nabla \phi_\bfv}{\phi_\bfv}
  &= \sum_{f \in F_\bfv} \bfp_f
     - \sum_{\bfu \in V} \phi_\bfu \sum_{f \in F_\bfu} \bfp_f \cr
  &= \sum_{f \in F_\bfv} \bfp_f
     - \sum_{f \in F} \mu_f \bfp_f \cr
  &= \sum_{f \in F_\bfv} (1-\mu_f) \bfp_f
  - \sum_{f \in F \setminus F_\bfv} \mu_f \bfp_f,
\end{align*}
and therefore,
$$ |\nabla \phi_\bfv| \le A_\bfv + B_\bfv, $$
where
$$ A_\bfv =  \phi_\bfv
   \sum_{f \in F_\bfv} (1-\mu_f) \frac{1}{h_f}, \qquad
   B_\bfv = \phi_\bfv \sum_{f \in F \setminus F_\bfv}
  \mu_f \frac{1}{h_f}. $$
Switching from summation over vertices-then-faces to faces-then-vertices, we have
\begin{equation}\label{eq:Av}
   \sum_{\bfv \in V} A_\bfv = 
   \sum_{f \in F} \sum_{\bfv \in V_f}
    \phi_\bfv (1-\mu_f) \frac{1}{h_f}
   = \sum_{f \in F} \mu_f (1-\mu_f) \frac{1}{h_f} =: C,
\end{equation}
and
\begin{equation}\label{eq:Bv}
   \sum_{\bfv \in V} B_\bfv = 
   \sum_{f \in F} \sum_{\bfv \in V \setminus V_f}
    \phi_\bfv \mu_f \frac{1}{h_f}
   = \sum_{f \in F} (1-\mu_f) \mu_f \frac{1}{h_f} = C,
\end{equation}
as well, and thus $\lambda\leq 2C$.
By Lemma~\ref{lem:hstarineq},
\begin{equation}
\label{eq:bd-C}
C \leq \sum_{f \in F} \mu_f \frac{1}{h_\ast}
   = \sum_{\bfv \in V} |F_\bfv| \phi_\bfv \frac{1}{h_\ast}
    = \frac{d}{h_\ast}.
\end{equation}
Therefore $\Lambda\leq 2d/h_\ast$ as claimed.
\end{proof}\\

The result of Theorem~\ref{thm:bound-gen-d} can be viewed as an improvement and generalization of a bound on $|\nabla\phi_\bfv|$ for Wachspress coordinates on polygons  given by Gillette, Rand and Bajaj in \cite[Lemma 6]{GRB1}.
Using Propositions 4, 7, and 8 from~\cite{GRB1}, we can write their bound as
\[ |\nabla\phi_\bfv| \leq \frac{\pi^2}{2}\left(\frac {4}{(d_\ast)^4\sin(\beta_\ast/2)\cos(\beta^\ast/2)\sin\beta^\ast}\right)^{2\beta^\ast/(\pi-\beta^\ast)},\]
where the interior angle $\beta_\bfv$ at vertex $\bfv$ is assumed to lie in $[\beta_\ast,\beta^\ast]\subset (0,\pi)$ and the length of an edge of $P$ is at least $d_\ast$. 
We now characterize our bound in terms of these same geometric quality measures to further illustrate the simplification provided here.
\begin{corollary}\label{cor:bound2dalt}
Let $P$ be a strictly convex polygon with minimum and maximum interior angles $\beta_\ast$ and $\beta^\ast$, respectively, and minimum edge length $d_\ast$.
Then
\begin{equation}\label{eq:bound2dalt}
 |\nabla \phi_\bfv| \le \frac{4}{d_\ast (\sin\beta_\ast)(\sin\beta^\ast)}.
\end{equation}
\end{corollary}

\begin{proof}
Label the vertices of $P$ in a counterclockwise fashion and let $\beta_k := \angle \bfv_{k-1} \bfv_k \bfv_{k+1}$ be the interior angle of $P$ at $\bfv_k$.
Either $\beta_k<\pi/2$ and $\sin(\beta_k)>\sin\beta_\ast$ or $\beta_k\geq\pi/2$ and $\sin(\beta_k)\geq \sin\beta^\ast$.
Since $P$ is strictly convex, we have $\sin\beta_\ast, \sin\beta^\ast\in (0,1]$ and hence $\sin(\beta_k)\geq (\sin\beta_\ast)(\sin\beta^\ast)>0$, without any qualification on $\beta_k$.
Now, again by the convexity of $P$, we have that
\[ \min_{\ell \ne k,k+1} h_k(\bfv_\ell) = \min(h_k(\bfv_{k-1}),h_k(\bfv_{k+2})).\]
Observe that
\[h_k(\bfv_{k-1}) 
    = |\bfv_{k-1} - \bfv_k| \sin(\beta_k)
   \geq d_\ast (\sin\beta_\ast)(\sin\beta^\ast),\]
and, similarly, $h_k(\bfv_{k+2}) \geq d_\ast (\sin\beta_\ast)(\sin\beta^\ast)$.
It follows that $h_\ast \geq d_\ast (\sin\beta_\ast)(\sin\beta^\ast)$.
By Theorem~\ref{thm:bound-gen-d}, $|\nabla\phi_\bfv|\leq 4/h_\ast$ and the result follows.
\end{proof}

\subsection{Simplices}

When $P$ is a simplex, we can improve the upper bound, as the following lemma demonstrates.\\

\begin{lemma}
\label{lem:smplx-up-bd}
Let $P$ be a $d$-simplex in $\RR^d$.  Then 
\[\Lambda \leq \frac{d+1}{h_\ast},\]
with equality in the case that $P$ is regular.
\end{lemma}\\
\begin{proof}
Note that  $|V| = d+1$. 
For any $\bfv\in V$ and $\bfx\in P$, we have $\nabla \phi_\bfv(\bfx) = -\bfn_f / h_f(\bfv)$, where $f$ is the unique face opposite to $\bfv$.
It follows that
\[ \Lambda = \sup_{\bfx \in P} \sum_{\bfv \in V} \frac {1}{h_f(\bfv)}  \leq \frac{d+1}{h_\ast}, \]
with equality in the case of a regular $d$-simplex.
\end{proof}

\subsection{Hyper-rectangles}

\begin{figure}[h]
\centering
\includegraphics[width=.5\textwidth]{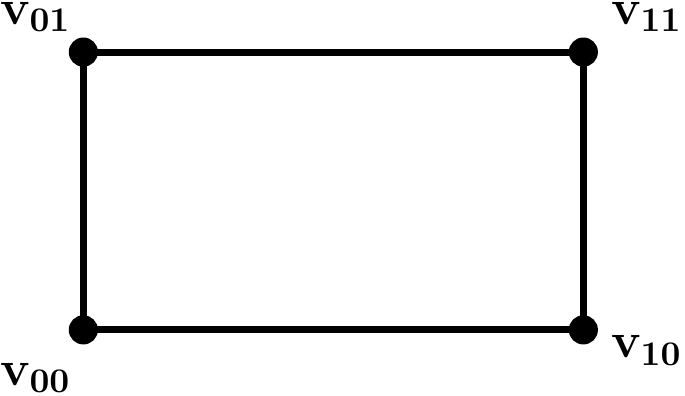}
\noindent \caption{Notation for the hyper-rectangle case in $d=2$: $\bfv_{00}=(a_1,a_2)$, $\bfv_{10}=(b_1,a_2)$, $\bfv_{01}=(a_1,b_2)$, $\bfv_{11}=(a_2,b_2)$.
}
\label{fg:sq-notation}
\end{figure}

We can also improve the upper bound in the case of a hyper-rectangle, i.e.\
\[ P = [a_1,b_1] \times [a_2,b_2] \times \cdots \times [a_d,b_d]\subset\RR^d. \]
In this case, the Wachspress coordinates are just the standard multilinear basis functions,  allowing for direct computation of their gradients and the following theorem.\\

\begin{theorem}
\label{theorem:hyperrect}
Let $P$ be a hyper-rectangle in $\RR^d$.  Then
\[ \Lambda \le \frac{\sqrt{d} + d}{h_\ast}, \]
with equality in the case that $P$ is a hypercube.
\end{theorem}\\

\begin{proof}
The $2^d$ vertices of $P$ can be indexed by $\eps \in \{0,1\}^d$, with $(\bfv_\eps)_i = a_i$ if $\eps_i = 0$ and  $(\bfv_\eps)_i = b_i$ if $\eps_i = 1$, $i=1,\ldots,d$.
The notation is shown in Figure~\ref{fg:sq-notation}.
Then, with $h_i := b_i - a_i$, we have
\[\phi_\eps(\bfx)
   = \prod_{i=1}^d
    \frac{1}{h_i} (x_i-a_i)^{\eps_i}
                (b_i-x_i)^{1-\eps_i}, \]
Since each factor is linear, we have that
\[ \frac{\partial}{\partial x_i}\phi_\eps(\bfx) = 
  (-1)^{1-\eps_i}\frac{1}{h_i} \prod_{j\not=i}  \frac{1}{h_j}(x_j-a_j)^{\eps_j}(b_j-x_j)^{1-\eps_j}
  ,\]
whence
\begin{equation}
\label{eq:hr-norm-exp}
|\nabla \phi_\eps(\bfx)|
   = \left(\prod_{i=1}^d
    \frac{1}{h_i}\right)\left(
     \sum_{i=1}^d
       \prod_{j \ne i}
       (x_j-a_j)^{2\eps_j}
                (b_j-x_j)^{2(1-\eps_j)}
     \right)^{1/2}.
\end{equation}
Now, observe that on an interval $[a,b]$ and for any $c,d \geq 0$, the two functions
\[ \sqrt{c (x-a)^2 + d}\quad\text{and}\quad \sqrt{c (b-x)^2 + d}\]
are convex for $x \in [a,b]$.
For any variable $x_i$, $|\nabla \phi_\eps(\bfx)|$ takes on one of these forms when the other $d-1$ variables are held fixed, meaning it is axially convex.
Since axial convexity is a closed property under addition, $\lambda$ is also axially convex.
We thus have
\[ \lambda(x_1,x_2,\ldots,x_d) \leq \max\left\{~
   \lambda(a_1,x_2,\ldots,x_d),
   \lambda(b_1,x_2,\ldots,x_d)~\right\}.\]
Similarly,
\begin{align*}
 \lambda(a_1,x_2,x_3,\ldots,x_d) & \leq \max\left\{~
   \lambda(a_1,a_2,x_3,\ldots,x_d),
   \lambda(a_1,b_2,x_3,\ldots,x_d)~\right\}, \\
 \lambda(b_1,x_2,x_3,\ldots,x_d) & \leq \max\left\{~
   \lambda(b_1,a_2,x_3,\ldots,x_d),
   \lambda(b_1,b_2,x_3,\ldots,x_d)~\right\}. 
\end{align*}
Thus, by induction, we arrive at
\[ \Lambda = \max_{\eps\in\{0,1\}^d} \lambda(\bfv_\eps). \]
Evaluating $\lambda(\bfv_\eps)$ using (\ref{eq:hr-norm-exp}), we see that
\[ |\nabla \phi_\eps(\bfv_\eps)|
   = \left( \sum_{i=1}^d \frac{1}{h_i^2} \right)^{1/2}. \]
The vertex $\bfv_\eps$ has $d$ neighbors, $\bfv_{\tau_i}$, where $\tau_i = (\eps_1,\ldots,1-\eps_i,\ldots,\eps_d)$.
We find that
\[ |\nabla \phi_{\tau_i}(\bfv_\eps)| 
= \frac{1}{h_i}. \]
Further, for fixed $\eps$, if $\eta \in \{0,1\}^d$ with
$\eta \not \in \{\eps,\tau_i\}$, then
$|\nabla \phi_\eta(\bfv_{\eps})| = 0$.
Therefore,
\[ \Lambda = \lambda(\bfv_{\epsilon}) =
   \left( \sum_{i=1}^d \frac{1}{h_i^2} \right)^{1/2} +~
          \sum_{i=1}^d \frac{1}{h_i}. \]
Since $h_\ast = \min_i h_i$, the desired inequality follows.
In the special case of a cube, all $h_i$ are equal, thereby completing the result.
\end{proof} \\

While Lemma~\ref{lem:smplx-up-bd} and Theorem~\ref{theorem:hyperrect} suggest that the upper bound from Theorem~\ref{thm:bound-gen-d} could possibly be improved, we will see in the next section that when $d=2$, the bound from the Theorem~\ref{thm:bound-gen-d} is in fact sharp over the class of simple convex polygons.

\section{Lower bounds}\label{sec:lower-bds}

Our lower bounds on $\Lambda$ are based on the observation that for any $\bfv\in V$, $\Lambda \geq \lambda(\bfv)$. 
We now broaden our scope from generalized Wachspress coordinates to allow any generalized barycentric coordinates satisfying the linear completeness property (\ref{eq:linprec3d}) and that are $C^1$ at the vertices of $\bfv$.  To be clear, we do \textit{not} assume that piecewise interpolation using the $\phi_\bfv$ is $C^1$, only that on a given polytope $P$, the associated functions $\phi_\bfv$ have well-defined gradients at the vertices of $P$.
We start with the polygonal case in order to clarify the subsequent generalizations to $d>2$.

\subsection{Polygons}

We start with a general lemma. 
Suppose $f:\RR^2 \to \RR$
is $C^1$ in a neighborhood of a point $\bfx \in \RR^2$.
Then $f$ has a directional derivative $D_\bfv f(\bfx)$ for any
non-zero vector $\bfv = (v_1,v_2) \in \RR^2$, which can be expressed
in terms of its gradient as
$$ D_\bfv f(\bfx) = \bfv \cdot \nabla f(\bfx). $$

\begin{lemma}
\label{lem:gradf-form}
If $\bfv$ and $\bfw$ are two linearly independent
vectors in $\RR^2$ then
\begin{equation}\label{eq:vw1}
 \nabla f
   = \frac{(D_\bfw f) \bfv^\perp - (D_\bfv f) \bfw^\perp}{\bfv \times \bfw},
\end{equation}
and
\begin{equation}\label{eq:vw2}
 |\nabla f| = 
     \frac{|(D_\bfw f) \bfv - (D_\bfv f) \bfw|}{|\bfv \times \bfw|},
\end{equation}
where $\bfv^\perp := (-v_2,v_1)$,
$\bfv \times \bfw := v_1 w_2 - v_2 w_1$, and
$|\cdot|$ is the Euclidean norm.
\end{lemma}\\

\begin{proof}
To show (\ref{eq:vw1}) it is sufficient to show that
$\bfv \cdot \nabla f = D_\bfv f$ and $\bfw\cdot\nabla f = D_\bfw f$.
These follow from taking the scalar product of the first
equation with $\bfv$ and $\bfw$ respectively, using the
fact that $\bfv \cdot \bfw^\perp = -\bfv \times \bfw$.
Equation (\ref{eq:vw2}) follows from the fact that $(D_\bfw f) \bfv^\perp - (D_\bfv f) \bfw^\perp = ((D_\bfw f) \bfv - (D_\bfv f) \bfw)^\perp$.
\end{proof}

We now use Lemma~\ref{lem:gradf-form} to estimate $\Lambda$ from below. 
Suppose that $P \in \RR^2$ is a convex polygon with vertices
$\bfv_1,\ldots,\bfv_n$ indexed in some counterclockwise ordering, and that $\phi_i:P \to \RR$, $i=1,\ldots,n$
are any set of generalized barycentric coordinates. 
Let $u:P\rightarrow\RR$ with $u_i:=u(\bfv_i)$ so that the formula (\ref{eq:I-def}) for $I(u)$ gives
\[ I(u) = \sum_{i=1}^n \phi_i(\bfx) u_i. \]
Observe that $I(u): \overbar P\rightarrow\RR$ is piecewise linear on $\partial P$.
Thus, at any vertex $\bfv_i$, we have 
\[ D_{\bfe_{i-1}}u(\bfv_i) = d_{i-1}\quad\text{and}\quad D_{\bfe_i}u(\bfv_i) = d_i\]
where $\bfe_j = \bfv_{j+1}-\bfv_j$ and $d_j = u_{j+1}-u_j$.
Letting $g:= I(u)$ to ease notation and evaluating (\ref{eq:vw2}) with $f=g$ at $\bfx=\bfv_i$ yields
\begin{equation}
\label{eq:nab-g-2d}
|\nabla g(\bfv_i)| = \frac{|d_i \bfe_{i-1} - d_{i-1}\bfe_i|} {\bfe_{i-1} \times \bfe_i}.
\end{equation}
\begin{corollary}
\label{cor:lam-bfv-form}
Let $P$ be a convex polygon as above.  Then
\[ \lambda(\bfv_i) =
   \frac{|\bfe_i| + |\bfe_i + \bfe_{i-1}| + |\bfe_{i-1}|}
                 {\bfe_{i-1} \times \bfe_i}. \]
\end{corollary}

\begin{proof}
Note that if $u=\phi_j$ then $u_i=\delta_{ij}$ so that $I(u)=\phi_j$. 
Using this fact and (\ref{eq:nab-g-2d}), we have the following formulae.
For $u = \phi_{i-1}$, $d_{i-1} = -1$, $d_i = 0$, and so
$$ |\nabla \phi_{i-1}(\bfv_i)| = 
   \frac{|\bfe_i|}
                 {\bfe_{i-1} \times \bfe_i}. $$
For $u = \phi_i$, $d_{i-1} = 1$, $d_i = -1$, and so
$$ |\nabla \phi_{i}(\bfv_i)| = 
   \frac{|\bfe_i + \bfe_{i-1}|}
                 {\bfe_{i-1} \times \bfe_i}. $$
For $u = \phi_{i+1}$, $d_{i-1} = 0$, $d_i = 1$, and so
$$ |\nabla \phi_{i+1}(\bfv_i)| = 
   \frac{|\bfe_{i-1}|}
                 {\bfe_{i-1} \times \bfe_i}. $$
Since $\nabla\phi_j(\bfv_i)=0$ for all other $j$, the result follows.
\end{proof}
\vspace{2mm}

\begin{theorem}
Let $P$ be the regular $n$-gon with vertices on the unit circle.  Then
\[\Lambda \ge \frac{2(1+ \cos(\pi/n))}{h_\ast} \qquad
\text{and}
\qquad
h_\ast = 4 \sin^2(\pi/n) \cos(\pi/n).\]
\end{theorem}\\

\emph{Proof.}
Let $\bfv_i = (\cos(i\theta),\sin(i\theta))$, $i=1,\ldots,n$, with $\theta = 2\pi/n$.
Since $\Lambda \ge \lambda(\bfv_n)$, it suffices to show that $\lambda(\bfv_n)$ is equal to the desired lower bound.
Since
\[ \bfv_{n-1} = (\cos\theta,-\sin\theta),
   \qquad \bfv_n = (1,0), \quad \bfv_1 = (\cos\theta,\sin\theta),\]
we have 
\[ \bfe_{n-1} = \bfv_n - \bfv_{n-1} = (1-\cos\theta,\sin\theta)\quad\text{and}\quad \bfe_n = \bfv_1 - \bfv_n = (\cos\theta-1,\sin\theta). \]
Therefore,
\begin{align*}
|\bfe_n| = |\bfe_{n-1}| & = 2\sin(\theta/2), \\
|\bfe_n + \bfe_{n-1}| = |\bfv_1 - \bfv_{n-1}| & = 2\sin(\theta),
\end{align*}
and
$$ \bfe_{n-1} \times \bfe_n = 2(1-\cos\theta)\sin\theta. $$
Thus by Corollary~\ref{cor:lam-bfv-form} and the double angle formulas,
$$ \lambda(\bfv_n) = \frac{2\sin(\theta/2) + \sin\theta}
                       {(1-\cos\theta)\sin\theta}
                   = \frac{1+ \cos(\theta/2)}
                       {2\sin^2(\theta/2) \cos(\theta/2)}. $$
Thus, it only remains to show that $h_\ast$ has the desired expression.
We have that
$$ h_\ast = h_1(\bfv_n) = (\bfv_1 - \bfv_n) \cdot \bfn_1. $$
Since $\bfn_1$ is the unit normal to the edge between $\bfv_1$ and $\bfv_2$, we have
$$ \bfn_1 = \begin{bmatrix} \cos(3\theta/2) \\[1mm] \sin(3\theta/2)\end{bmatrix}. $$
Using summation and sum-to-product trigonometric  identities, we compute that
$$ h_\ast = (\cos\theta-1) \cos(3\theta/2) 
    + \sin\theta \sin(3\theta/2)
  = 4 \sin^2(\theta/2) \cos(\theta/2). \qquad \endproof $$

We obtain a further lower bound on $\Lambda$, this time for general convex polygons.\\

\begin{theorem}
\label{thm:lwr-bd-polyg}
Let $P$ be a convex polygon in $\RR^2$.  For any generalized barycentric coordinates  $\phi_1,\ldots,\phi_n$ on $P$ that are $C^1$ at the vertices of $P$, 
\[ \Lambda \ge \frac{1}{h_\ast}. \]
Further, if all the interior angles of $P$ are obtuse,
\[ \Lambda \ge \frac{2}{h_\ast}. \]
\end{theorem}\\

\begin{proof}
There is some index $i$ such that either $h_\ast = h_i(\bfv_{i-1})$ or $h_\ast = h_{i-1}(\bfv_{i+1})$.
Without loss of generality, assume that $h_\ast = h_i(\bfv_{i-1})$.
Then, with $\beta_i$ the interior angle of $P$ at $\bfv_i$,
whether acute or obtuse, $h_i(\bfv_{i-1}) = |\bfe_{i-1}| \sin\beta_i,$ and so
\[ \bfe_{i-1} \times \bfe_i = |\bfe_{i-1}| |\bfe_i| \sin\beta_i  = |\bfe_i| h_\ast.\]
Then by Corollary~\ref{cor:lam-bfv-form}, $ \lambda(\bfv_i) \ge 1/h_\ast.$
In the case that $\beta_i$ is obtuse, i.e., $\beta_i \ge \pi/2$, $ |\bfe_{i-1} + \bfe_i| = |\bfv_{i+1}-\bfv_{i-1}| \ge |\bfe_i|,$ and Corollary~\ref{cor:lam-bfv-form}  gives $ \lambda(\bfv_i) \ge 2/h_\ast.$
\end{proof}

\subsection{Simple polyhedra}

We now consider the case where $P$ is a simple convex polyhedron in $\RR^3$ and $\phi_\bfv$ are any set of generalized barycentric coordinates on $P$ that are $C^1$ at the vertices of $P$.
We first prove a lemma about the geometry of $P$.\\

\begin{lemma}\label{lem:hstarattained}
For any face $f \in F$, let $\bfw \in V \setminus V_f$ be a vertex of $P$ satisfying
\begin{equation}\label{eq:wmin}
 h_f(\bfw) = \min_{\bfu \in V \setminus V_f} h_f(\bfu).
\end{equation}
Then $\bfw$ is a neighbor of a vertex in $f$.
\end{lemma}\\
\begin{proof}
We show first that any vertex $\bfu  \in V \setminus V_f$
has a neighbor $\bfv \in N_\bfu$ such that $h_f(\bfv) < h_f(\bfu)$.
Indeed, if to the contrary,
$h_f(\bfv) \ge h_f(\bfu)$ for all $\bfv \in N_\bfu$
then, by the convexity of $P$,
the half-space $H \subset \RR^3$ defined by the equation
$$ \{ \bfy  \in \RR^3 : h_f(\bfy) \ge h_f(\bfu) \} $$
must contain $P$, which contradicts the fact that the
vertices of $f$ lie outside $H$.

Suppose now that $\bfw \in V \setminus V_f$ is any vertex satisfying
(\ref{eq:wmin}).
By the above argument, there must be some neighbor $\bfv \in N_\bfw$
such that $h_f(\bfv) < h_f(\bfw)$.
If $\bfv \not\in V_f$,  this contradicts the definition of $\bfw$.
Therefore $\bfv \in V_f$, which proves the lemma.
\end{proof}\\

Similar to the planar case, if $\bfu$, $\bfv$, $\bfw$ are
linearly independent vectors in $\RR^3$,
then for a smooth enough function $f:\RR^3\to \RR$,
\begin{equation}\label{eq:vw13D}
 \nabla f
   = \frac{(D_\bfw f) \bfu\times \bfv
            + (D_\bfu f) \bfv\times \bfw
            + (D_\bfv f) \bfw\times \bfu}
            {\det(\bfu,\bfv,\bfw)}.
\end{equation}

Let $u:\overbar P\rightarrow\RR$ with $u_\bfv:=u(\bfv)$ so that the formula (\ref{eq:I-def}) for $I(u)$ gives
\[ I(u)(\bfx) = \sum_{\bfv\in V} \phi_\bfv(\bfx) u_\bfv. \]
Fix $\bfv \in V$ and let $\bfw_i$, $i=1,2,3$ be the three neighbors of $\bfv$, in some clockwise ordering as seen from outside $P$.
Let $\bfe_i = \bfw_i - \bfv$ and $d_i = u_{\bfw_i} - u_\bfv$.
Since $I(u)$ is linear along the edges of $P$, $D_{\bfw_i}u(\bfv) = d_i$.
Letting $g:= I(u)$ to ease notation and evaluating (\ref{eq:vw13D}) with $f=g$ at $\bfx=\bfv$ yields the formula
\begin{equation}
\label{eq:nab-g-3d}
 \nabla g(\bfv)
   = \frac{d_1 \bfe_2\times \bfe_3
            + d_2\bfe_3\times \bfe_1
            + d_3\bfe_1\times \bfe_2}
            {\det(\bfe_1,\bfe_2,\bfe_3)}.
\end{equation}
We use this to prove a lower bound on $\Lambda$ for convex polyhedra.\\

\begin{theorem}
\label{thm:lwr-bd-polyhedr}
Let $P$ be a convex polyhedron in $\RR^3$.  For any generalized barycentric coordinates  $\phi_\bfv$ on $P$ that are $C^1$ at the vertices of $P$, 
\[ \Lambda \ge \frac{1}{h_\ast}. \]
\end{theorem}
\begin{proof}
By Lemma~\ref{lem:hstarattained}, we can set $f\in F$ and $\bfv\in f$ so that $h_\ast=h_f(\bfw_3)$. 
Since $\phi_\bfu(\bfv) = \delta_{\bfu\bfv}$, 
\[ \lambda(\bfv) = |\nabla \phi_\bfv(\bfv)|
             + \sum_{i=1}^3 |\nabla \phi_{\bfw_i}(\bfv)|. \]
By (\ref{eq:nab-g-3d}), 
\[
 \nabla \phi_\bfv(\bfv)
   = \frac{-(\bfe_2\times \bfe_3
            + \bfe_3\times \bfe_1
            + \bfe_1\times \bfe_2)}
            {\det(\bfe_1,\bfe_2,\bfe_3)},
\]
and
\[
 \nabla \phi_{\bfw_1}(\bfv)
   = \frac{\bfe_2\times \bfe_3}
            {\det(\bfe_1,\bfe_2,\bfe_3)}, \quad
 \nabla \phi_{\bfw_2}(\bfv)
   = \frac{\bfe_3\times \bfe_1}
            {\det(\bfe_1,\bfe_2,\bfe_3)}, \quad
 \nabla \phi_{\bfw_3}(\bfv)
   = \frac{\bfe_1\times \bfe_2}
            {\det(\bfe_1,\bfe_2,\bfe_3)}.
\]
Note that $f$ is the face of $P$ containing the
vertices $\bfv$, $\bfw_1$, and $\bfw_2$.
We thus have
\[ \det(\bfe_1,\bfe_2,\bfe_3)
   = (\bfe_1\times \bfe_2) \cdot \bfe_3
   = |\bfe_1\times \bfe_2| (-\bfn_f) \cdot \bfe_3,\]
Using this and (\ref{eq:def-hf}), we see that  $|\nabla \phi_{\bfw_3}(\bfv)| = 1 / {|\bfe_e\cdot\bfn_f|} = {1} / {h_f(\bfw_3)}$, and hence
\[ \lambda(\bfv) \ge \frac{1}{h_f(\bfw_3)} = \frac 1{h_\ast}. \]
Since $\Lambda\geq\lambda(\bfv)$ for any $\bfv\in V$, we have completed the proof.
\end{proof}

\subsection{Simple polytopes}

For a simple convex polytope $P$ in $\RR^d$ with $d>3$, the bound $\Lambda\geq 1/h_\ast$ holds by the same analysis as in the polyhedral case just presented.
Note that the proof technique of Lemma~\ref{lem:hstarattained} is not specific to $\RR^3$ and thus carries over to the generic $d$ case immediately.
The proof technique for Theorem~\ref{thm:lwr-bd-polyhedr} is also not specific to $\RR^3$, although the notation becomes more dense.
Since the cases $d>3$ are of less interest from an application perspective, and in the interest of space, we omit stating the straightforward generalizations of the proof.

\section{Numerical Experiments}\label{sec:results}

As discussed in the introduction, the upper bounds derived on $\Lambda$ suffice to ensure that a Lagrange-style finite element method employing generalized Wachspress coordinates as basis functions will obtain a linear order \textit{a priori} error estimate as stated in (\ref{eq:opt-conv-est}).
We provide numerical evidence of this convergence not only to confirm its theoretical validity, but also to demonstrate that implementing such basis functions is computationally viable. 
The numerical computations were carried out 
in \texttt{MATLAB\texttrademark}.  
We include code for the computation of the basis functions on polygons and polyhedra in the Appendix.

Consider the following three-dimensional Poisson boundary-value problem:
\begin{subequations}\label{eq:poisson}
\begin{align}
-\nabla^2 u &= f \qquad \textrm{in } \Omega = (0,1)^3, \\
\label{eq:ebc}
u &= 0 \qquad \textrm{on } \partial \Omega.
\end{align}
\end{subequations}
The weak form of this problem is: find $u \in H^1_0(\Omega)$ such that
\begin{align}
\int_\Omega \nabla u \cdot \nabla w \, d\bfx = \int_\Omega f w \, d\bfx,\quad \forall w\in H_0^1(\Omega), \label{eq:weak}
\end{align}
where $H_0^1(\Omega)$ is the standard
Sobolev space of degree $1$ with vanishing values on the boundary.

\begin{figure}
\centering
\begin{tabular}{cccc}
  \hspace{-1.0cm}\includegraphics[width=4.5cm,clip=]{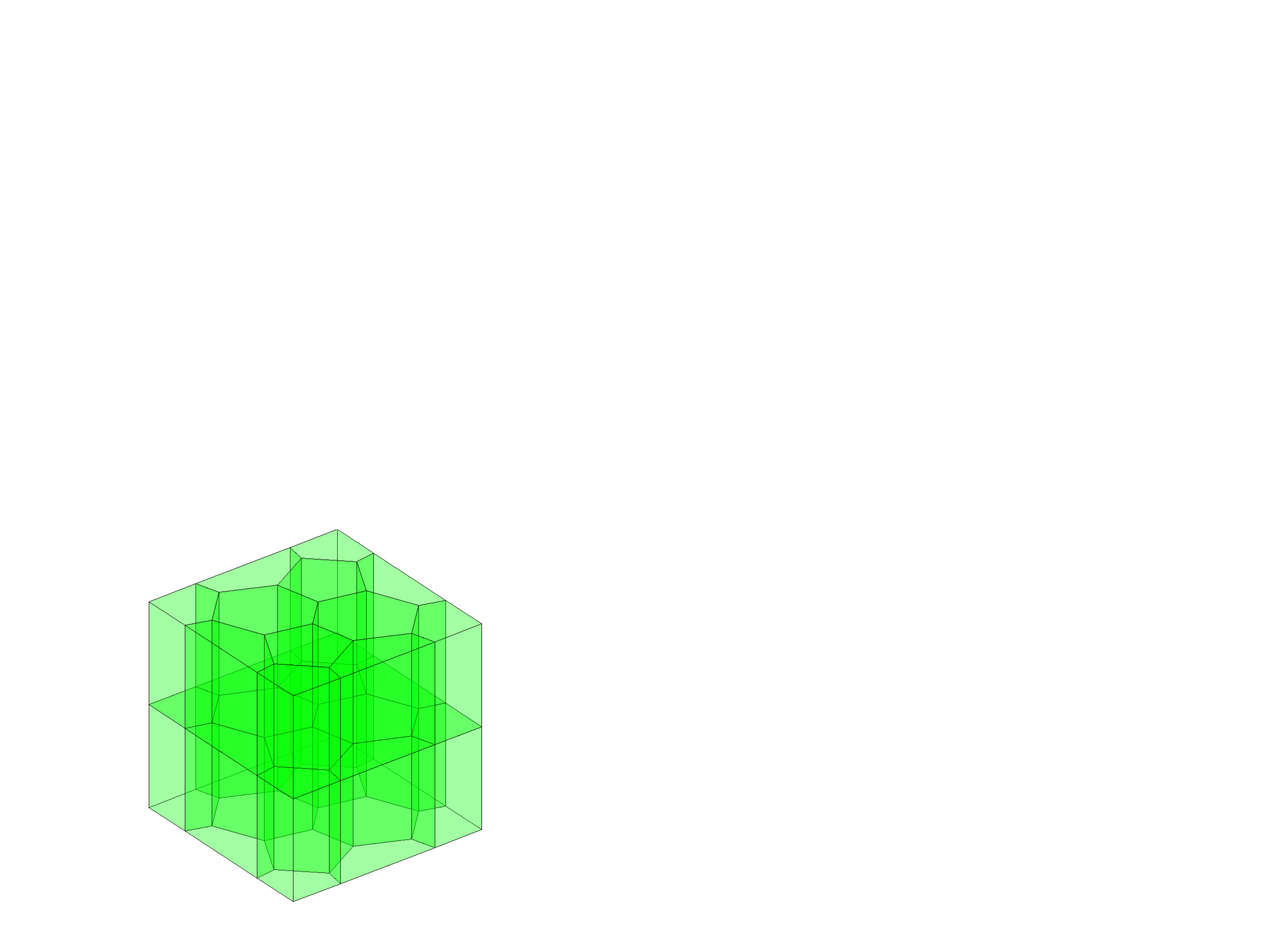} &
  \hspace{-1.4cm}\includegraphics[width=4.5cm,clip=]{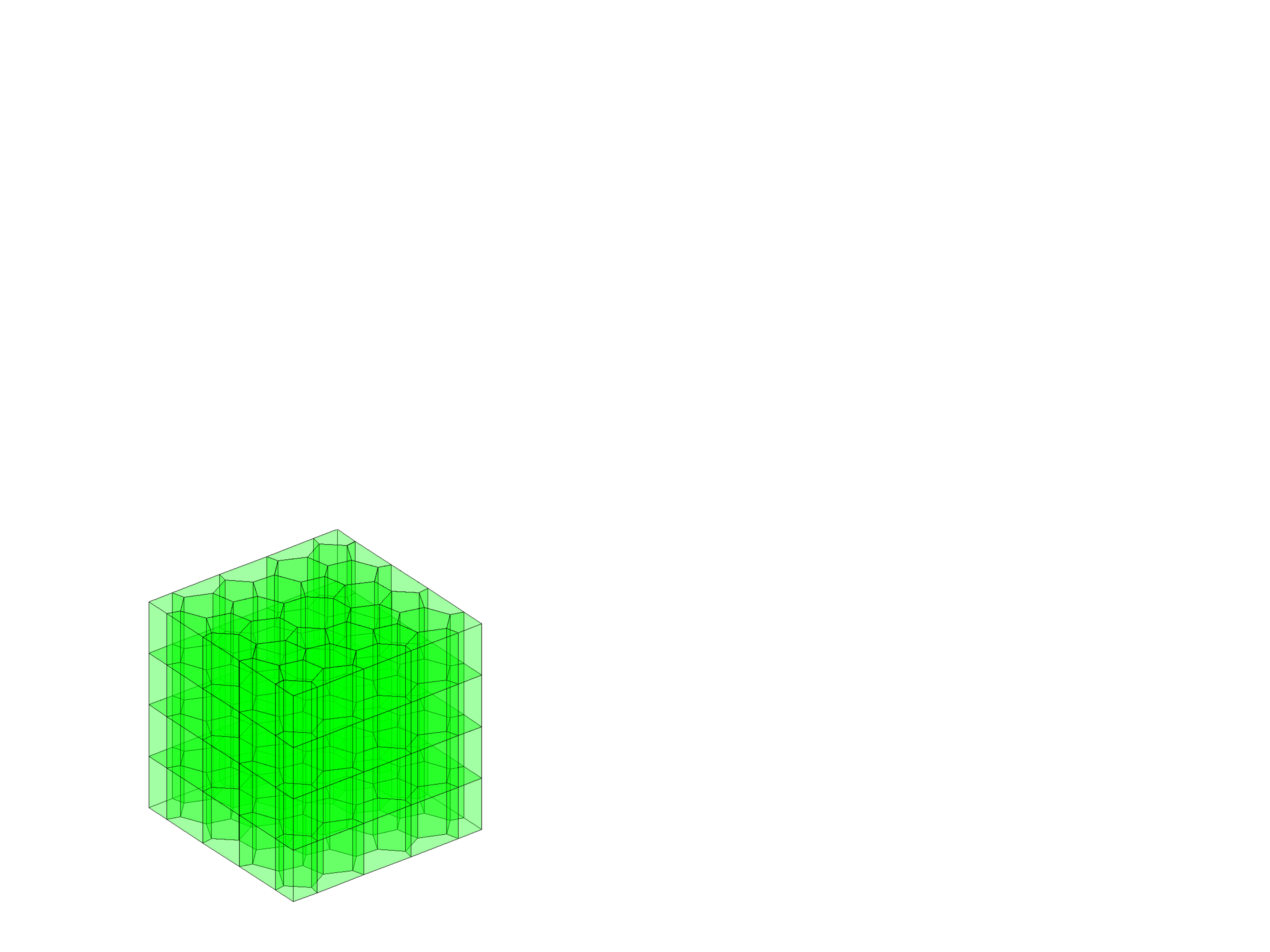} &
  \hspace{-1.4cm}\includegraphics[width=4.5cm,clip=]{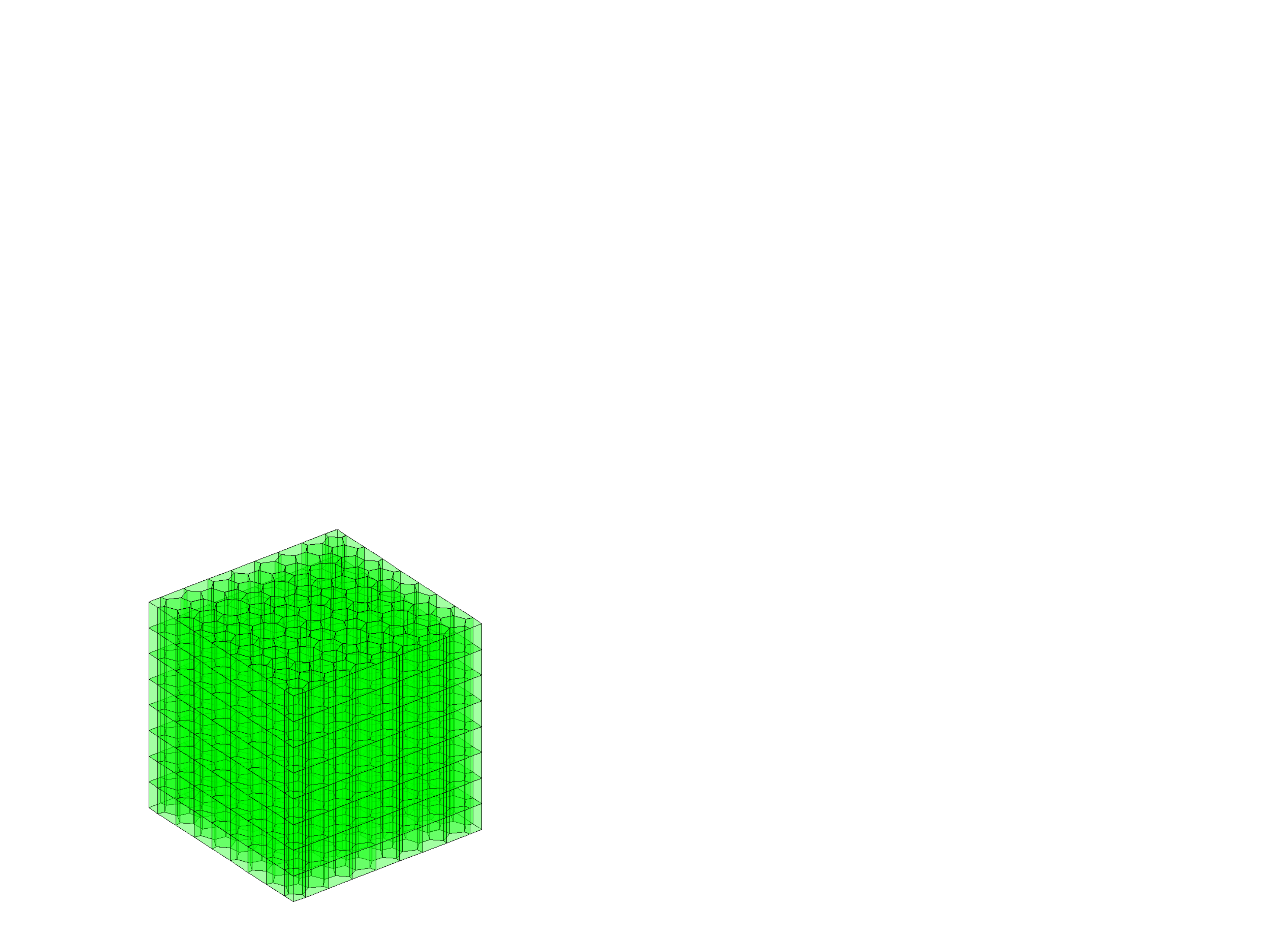} &
  \hspace{-1.4cm}\includegraphics[width=4.5cm,clip=]{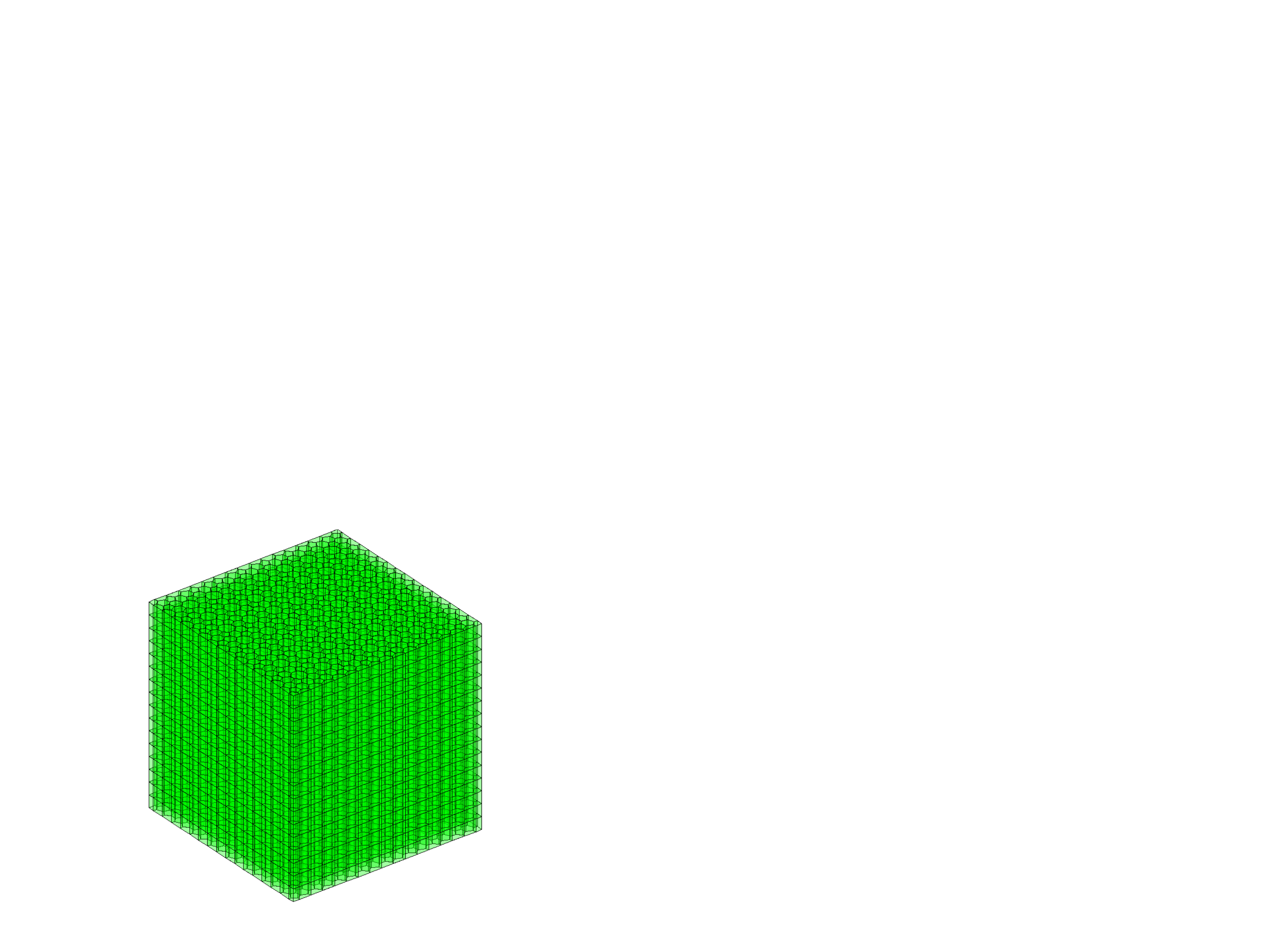} \\
%
  \hspace{-1.0cm} (a) &
  \hspace{-1.4cm} (b) &
  \hspace{-1.4cm} (c) &
  \hspace{-1.4cm} (d) \\ 
\end{tabular}
\caption{Polyhedral meshes on which we solve a Poisson problem using generalized Wachspress basis functions in a finite element method. A fifth mesh (e) of even finer resolution is not shown. }
\label{fig:meshes}
\end{figure}

\begin{table}
\caption{\noindent Relative error norms for solution of the Poisson problem on the polyhedral meshes shown in Figure~\ref{fig:meshes}.  The quantity $h$ denotes the maximum diameter of a mesh element.  }\label{tab:conv}
\centering
\vspace*{0.1in}
\begin{tabular}{ | c | c | c | c | c | c | c | }
\hline
Mesh & \# of nodes & $h$ & $\dfrac{||u - u^h||_{0,P}}{||u||_{0,P}}$ & Rate & $\dfrac{|u - u^h|_{1,P}}{|u|_{1,P}}$ & Rate \\ 
 \hline
a & 78     & 0.7071 &  $2.0 \times 10^{-1}$ & --  & $4.1 \times 10^{-1}$ & -- \\ \hline 
b & 380    & 0.3955 &  $5.4 \times 10^{-2}$ & 2.28 & $2.1 \times 10^{-1}$ & 1.14 \\ \hline 
c & 2340   & 0.1977 &  $1.4 \times 10^{-2}$ & 1.96 & $1.1 \times 10^{-1}$ &0.97 \\ \hline 
d & 16388  & 0.0989 &  $3.5 \times 10^{-3}$ & 1.99 & $5.4 \times 10^{-2}$ &0.99 \\ \hline 
e & 122628 & 0.0494 &  $8.8 \times 10^{-4}$ & 2.00 & $2.7 \times 10^{-2}$ &0.99 \\\hline 
\end{tabular}
\end{table}

We define a Galerkin method for approximating the solution to \eref{eq:weak}.
First, we fix a sequence of meshes made of polygonal prismatic elements, as shown in Figure~\ref{fig:meshes}.
The meshes have increasingly smaller values of $h$, the maximum diameter of a mesh element, but the geometric quality of the elements does not degrade with $h$.
More precisely, there is a constant $\gamma>0$ such that if any element from any mesh in the sequence is scaled to have diameter 1, the computed value of $h_\ast$ will be $\geq\gamma$.

For each simple polyhedral element $P_e$, we label the vertices $\bfv_1,\ldots,\bfv_n$, and associate the basis function $\phi_i(\vx) := \phi_\bfv(\vx)$ to the vertex $\bfv_i$.  
Note that different elements need not have the same number of vertices.
We choose the set of generalized Wachspress coordinates $\{\phi_i (\vx)\}_{i=1}^n$ to be the basis for the trial and test spaces on $P_e$.   This yields the linear system
\begin{equation}\label{eq:kefe}
\vm{K}_e \vm{d}_e = \vm{f}_e, \quad\text{with}\quad
\vm{K}_{e}^{i j} = \int_{P_e} \nabla \phi_i \cdot \nabla \phi_j \, 
d\vx , \quad
\vm{f}_{e}^{i} = \int_{P_e} f \phi_i \, d\vx ,
\end{equation}
and $\vm{d}_e$ is a vector of unknown coefficients. 
The global stiffness matrix $\vm{K}$ is formed by assembling contributions from all the $\vm{K}_e$; similarly, the global element source vector $\vm{f}$ is formed by assembling  contributions from all the $\vm{f}_e$. 
The global linear system, $\vm{K} \vm{d} = \vm{f}$, is solved after imposing the homogeneous Dirichlet boundary conditions.

We choose $f(\vx)$ in \eref{eq:poisson} so that the exact solution to the Poisson problem is $u(\vx) = xyz(1-x)(1-y)(1-z)$. 
To evaluate $\vm{K}_e$ and $\vm{f}_e$ in~\eref{eq:kefe}, we partition each polyhedral element into tetrahedra, and use a second-order accurate polynomial-precision quadrature rule (four quadrature points) within each tetrahedron~\cite{Shunn2012}. 
The relative error norms are listed in~\tref{tab:conv}.
The polyhedral finite element method delivers optimal convergence rates of $2$ and $1$ in the $L^2$ norm and the $H^1$ seminorm, respectively, as expected from the \textit{a priori} estimate (\ref{eq:opt-conv-est}).
 
\section{Conclusion and Future Directions}\label{sec:conc}

The results presented in this paper help to answer a significant question in finite element theory: what makes a good linear finite element?
While the notion of a `good' element is highly dependent on the application context, it is generally acknowledged that avoiding large angles~\cite{BA1976,J1976}  or large circumradii~\cite{K1991,R2012} is desirable for finite element methods.
The circumradius $r_\text{circ}$ of a triangle is defined as the radius of its circumcircle.
Using results from Shewchuk~\cite{S2002} for geometric properties of triangles, we find that $r_\text{circ}=\ell_\text{min} \ell_\text{med}/2h_\ast$ where $\ell_\text{min}$, $\ell_\text{med}$ are the lengths of the shortest two sides of the triangle.
Scaling the largest edge of the triangle to length 1, we see that
\[ r_\text{circ} = \frac{\ell_\text{min} \ell_\text{med}}{2h_\ast}\leq \frac{(\ell_\text{min}+\ell_\text{med})^2}{8h_\ast}\leq \frac 1{2 h_\ast}.\]
Therefore, a triangular mesh that avoids small $h_\ast$ values also avoids large circumradii.

The question of what makes a `good' linear finite element  becomes much more difficult for polytopes as their geometry can be quite exotic and difficult to characterize by only a few parameters.
Even on tetrahedra, this question is quite subtle as evidenced by the wealth of literature in this area.
Hence, our work is only a first viewpoint toward a generic characterization of the relationship between polytope element geometry and interpolation error estimates.
The quantity $h_\ast$ is easily computed and relates closely to the operator norm of the interpolant $I$, as evidenced both by our bounds on $\Lambda$ and by our numerical experiments.

Additionally, as discussed at the beginning of Section~\ref{sec:lower-bds}, recall that our lower bounds on $\Lambda$ pertain to any set of generalized barycentric coordinates that are $C^1$ at the vertices of simple polytopes.
Our proof technique reveals that no generalized barycentric coordinates can be $C^1$ at a non-simple vertex (i.e.\ a vertex incident to more than $d$ faces of $P$) as this would require the gradient at the vertex to satisfy more than $d$ linearly independent constraints.
Further, while the Wachspress coordinates are indeed $C^1$ at simple vertices, other coordinates may not be, leaving open the possibility that $\Lambda$ may in fact be smaller than $1/h_\ast$ for a different choice of generalized barycentric coordinates~\cite{F2003,FHK,FKR,JSW,S2004}.
We plan to pursue this question in future work.

Finally, the code used in our experiments, which appears in the Appendix, is designed to simplify the process of computing and integrating the Wachspress coordinate functions and their gradients on convex polygons and generic (i.e.\ simple or non-simple) convex polyhedra.
The code follows the definitions and derivations of Section~\ref{sec:upper-bds}.
We were motivated to develop this code after learning at a recent workshop~\cite{GBC2012} of the growing interest in the computer graphics and finite element communities for efficient implementation of polyhedral generalized barycentric coordinates.

\newpage

\section*{Acknowledgments}
The authors thank Gianmarco Manzini for providing the polyhedral meshes that
were used in the convergence study and Alexander Rand for helpful conversations 
about the error estimation literature.
AG acknowledges support as a postdoc from NSF Award DMS-0715146 and by NBCR at UC San Diego.
NS acknowledges support from NSF Award CMMI-1334783 at UC Davis.\\

\noindent
The published version of this paper will appear in SIAM Journal on Numerical Analysis in 2014.


%

\section*{Appendix A}
\texttt{MATLAB\texttrademark} codes for the computation of Wachspress basis functions on
convex polygons and polyhedra are listed. For geometric computations on 
polyhedra, the \texttt{geom3d} library is used~\cite{geom3d}.

\subsection*{Wachspress basis functions on convex polygons}

We restate the formulae for the Wachspress coordinates from \eref{eq:lambda-def} and \eref{eq:wachsimple} and their gradients from \eref{eq:def-R} and \eref{eq:R3d} in the case $d=2$ before giving the code used to compute them.

Let $P \subset \RR^2$ be a convex polygon, with vertices $\bfv_1,\ldots,\bfv_n \in \RR^2$ in some counter-clockwise ordering.
Let $\bfn_i \in \RR^2$ be the outward unit normal to the edge $\bfe_i = [\bfv_i,\bfv_{i+1}]$, with vertices indexed cyclically, i.e., $\bfv_{n+1} := \bfv_1$ etc.
For any $\bfx$ in $P$, let $h_i(\bfx)$ be the perpendicular distance of $\bfx$ to the edge $\bfe_i$ and define $\bfp_i(\bfx) := \bfn_i / h_i(\bfx)$.
Then the coordinate functions $\phi_i := \phi_{\bfv_i} : P \to \RR$ are given by
\begin{equation}\label{eq:phi2d}
\phi_i = \frac{ w_i }{\displaystyle \sum_{j=1}^n w_j }
\quad
\text{ where }
\quad
 w_i := \det(\bfp_{i-1}, \bfp_{i}).
\end{equation}
The gradient functions are given by
\begin{equation}\label{eq:nabphi2d}
\nabla \phi_i = \phi_i (\bfR_i - \sum_{j=1}^n \phi_j \bfR_j)
\quad
\text{ where }
\quad
\bfR_i := \bfp_{i-1} + \bfp_i.
\end{equation}
These functions are implemented in  \texttt{MATLAB\texttrademark} as follows.
\vspace{.5cm}

\begin{verbatim}
function [phi dphi] = wachspress2d(v,x)
%
% Evaluate Wachspress basis functions and their gradients in a convex polygon
%
% Inputs:
% v    : [x1 y1; x2 y2; ...; xn yn], the n vertices of the polygon in ccw
% x    : [x(1) x(2)], the point at which the basis functions are computed
% Outputs:
% phi  : output basis functions = [phi_1; ...; phi_n]
% dphi : output gradient of basis functions = [dphi_1; ...; dphi_n]

n = size(v,1);
w = zeros(n,1);
R = zeros(n,2);
phi  = zeros(n,1);
dphi = zeros(n,2);

un = getNormals(v);
\end{verbatim}
\newpage
\begin{verbatim}
p = zeros(n,2);
for i = 1:n
  h = dot(v(i,:) - x,un(i,:));
  p(i,:) = un(i,:) / h;
end

for i = 1:n
  im1 = mod(i-2,n) + 1;
  w(i) = det([p(im1,:);p(i,:)]);
  R(i,:) = p(im1,:) + p(i,:);
end

wsum = sum(w); 
phi = w/wsum;

phiR = phi' * R;
for k = 1:2
  dphi(:,k) = phi .* (R(:,k) - phiR(:,k));
end

function un = getNormals(v)
% Function to compute the outward unit normal to each edge

n = size(v,1);
un = zeros(n,2);
for i = 1:n
  d = v(mod(i,n)+1,:) - v(i,:);
  un(i,:) = [d(2) -d(1)]/norm(d);
end
\end{verbatim}

\subsection*{Wachspress basis functions on convex polyhedra}

We state the formulae for generalized Wachspress coordinates and their gradients on generic (simple or non-simple) convex polyhedra before giving the code used to compute them.
In the case of a simple polyhedron, this definition recovers the $d=3$ case of  \eref{eq:lambda-def}, \eref{eq:wachsimple}, \eref{eq:def-R} and \eref{eq:R3d}.

Let $P \subset \RR^3$ be a convex polyhedron, with vertex set $V$.
For each $\bfv\in V$, denote the $k \geq 3$ faces incident on $\bfv$ by $\bff_1, \bff_2, \ldots, \bff_k$, in some counter-clockwise order as seen from outside $P$.
Let $\bfn_1, \bfn_2, \ldots, \bfn_k$ be the outward unit normals to these faces, respectively.
Let $h_i(\bfx) > 0$ be the perpendicular distance of $\bfx$ from the face $\bff_i$ and define $\bfp_i(\bfx) := \bfn_i / h_i(\bfx)$.
Then the coordinate functions $\phi_{\bfv}: P \to \RR$ are given by 
\begin{equation}
\label{eq:phi3d}
\phi_\bfv = \frac{ w_\bfv }{\displaystyle \sum_{\bfu\in V} w_\bfu }
\quad
\text{ where }
\quad
w_\bfv :=  \sum_{i=1}^{k-2} w_{i,\bfv}
\quad
\text{ and }
\quad
w_{i,\bfv} := \det(\bfp_{i}, \bfp_{i+1}, \bfp_{k}).
\end{equation}
The gradient functions are then given by
\begin{equation}\label{eq:nabphi3d}
\nabla \phi_\bfv = \phi_\bfv \left(\bfR_\bfv - \sum_{\bfu\in V} \phi_\bfu \bfR_\bfu\right)
\quad
\text{ where }
\quad
\bfR_\bfv := \frac{1}{w_\bfv} \sum_{i=1}^{k-2} w_{i,\bfv} (\bfp_i + \bfp_{i+1} + \bfp_k).
\end{equation}
These functions are implemented in  \texttt{MATLAB\texttrademark} as follows.
\vspace{.5cm}

\begin{verbatim}
function [phi dphi] = wachspress3d(v,g,un,x)

% Evaluate Wachspress basis functions and their gradients
% in a convex polyhedron
%
% Inputs:
% v    : [x1 y1 z1; x2 y2 z2; ...; xn yn zn], the n vertices of the polyhedron
% g    : cell array: [i1 i2 ... i_{k1}]; . . . ; [i1 i2 ... i_{kn}],
%        which are the n neighborhood graphs, each in some counter-clockwise order
%        as seen from outside the polyhedron
% un   : [x1 y1 z1; x2 y2 z2; ...; xm ym zm], unit normal to each facet
% x    : [x(1) x(2) x(3)], the point at which the basis functions are computed
% Outputs:
% phi  : basis functions = [phi_1; ...; phi_n]
% dphi : gradient of basis functions = [dphi_1; ...; dphi_n]

n = size(v,1);
w = zeros(n,1);
R = zeros(n,3);
phi  = zeros(n,1);
dphi = zeros(n,3);

for i = 1:n
  f = g{i};
  k = length(f);
  p = zeros(k,3);
  for j = 1:k
    h = dot(v(i,:) - x,un(f(j),:));
    p(j,:) = un(f(j),:) / h;
  end

  wloc = zeros(k-2,1);
  Rloc = zeros(k-2,3);
  for j = 1:k-2
    wloc(j) = det([p(j,:); p(j+1,:); p(k,:)]);
    Rloc(j,:) = p(j,:) + p(j+1,:) + p(k,:);
  end

  w(i) = sum(wloc);
  R(i,:) = (wloc' * Rloc) / w(i);
end

wsum = sum(w); 
phi  = w/wsum;

phiR = phi' * R;
for d = 1:3
  dphi(:,d) = phi .* (R(:,d) - phiR(:,d));
end
\end{verbatim}


\end{document}